\documentclass[12pt]{article}
\usepackage{amsmath,amssymb,amsfonts}
\newtheorem{theorem}{Theorem}
\newtheorem{lemma}{Lemma}

\date{}
\title{\bf A Note on Goldbach Partitions of Large Even Integers}

\author{{Ljuben Mutafchiev}\\
\small American University in Bulgaria, 2700 Blagoevgrad, Bulgaria
\\[-0.8ex]
\small and Institute of Mathematics and Informatics of the
\\[-0.8ex]
\small Bulgarian Academy of Sciences \\[-0.8ex]
\small \tt {ljuben@aubg.bg}}
\date{\small Submitted: July 18, 2014; Accepted: Jan ..., 2015
}

\begin{document}
\maketitle

\begin{abstract}
Let $\Sigma_{2n}$ be the set of all partitions of the even
integers from the interval $(4,2n], n>2,$ into two odd prime
parts. We show that $\mid\Sigma_{2n}\mid\sim 2n^2/\log^2{n}$ as
$n\to\infty$. We also assume that a partition is selected
uniformly at random from the set $\Sigma_{2n}$. Let $2X_n\in
(4,2n]$ be the size of this partition. We prove a limit theorem
which establishes that $X_n/n$ converges weakly to the maximum of
two random variables which are independent copies of a uniformly
distributed random variable in the interval $(0,1)$. Our method of
proof is based on a classical Tauberian theorem due to Hardy,
Littlewood and Karamata. We also show that the same asymptotic
approach can be applied to partitions of integers into an
arbitrary and fixed number of odd prime parts.
\end{abstract}

\vspace{.2cm}

\section{Introduction and Statement of the Main Result}

For a given sequence of positive integers
$\Lambda=\{\lambda_1,\lambda_2,...\}$, by a $\Lambda$-partition of
the positive integer $n$, we mean a way of writing it as a sum of
positive integers from $\Lambda$ without regard to order; the
summands are called parts. Let $\mathcal {P}=\{p_1,p_2,...\}$ be
the sequence of all odd primes arranged in increasing order. A prime
partition is a $\Lambda$-partition with $\Lambda=\mathcal{P}$. Let
$Q(n)$ be the number of prime partitions of $n$. Hardy and
Ramanujan [6,7] were apparently the first who studied the
asymptotic behavior of the number of integer
($\Lambda=\{1,2,...\}$) and prime partitions for large $n$. For
prime partitions they proved the following asymptotic formula:
$$
\log{Q(n)}\sim 2\pi\sqrt{\frac{n}{3\log{n}}}, \quad n\to\infty.
 $$
The study of the asymptotic behavior of $Q(n)$ itself is quite
complicated. It turns out that the corresponding asymptotic
formula contains transcendental sums over the primes which can be
expressed in terms of zeros of the Riemann zeta function (for more
details, see e.g. [9; p. 240]). Recently Vaughan [16] proposed and
studied a modification of the problem, where $n$ is replaced by a
continuous real variable. His asymptotic results avoid
transcendental sums over primes.

Consider now the number $Q_m(n)$ of prime partitions of $n$ into
$m$ parts ($1\le m\le n$). The bivariate generating function of
the numbers $Q_m(n)$ is of Euler's type, namely,
\begin{equation}\label{euler}
G(x,z) =1+\sum_{n=1}^\infty z^n\sum_{m=1}^n Q_m(n)x^m
=\prod_{p_k\in\mathcal{P}}(1-xz^{p_k})^{-1}
\end{equation}
(the proof may be found in [1; Section 2.1]). In this note we
focus on the asymptotic behavior of the coefficients $Q_2(n)$ of
$x^2$ and $z^n$ in the power series expansion of $G(x,z)$ in
powers of $x$ and $z$. For $n>4$, $Q_2(n)$ counts the number of
ways of representing $n$ as a sum of two odd primes. Obviously,
$Q_2(n)=0$ if $n$ is odd. In 1742 Goldbach conjectured that
$Q_2(n)\ge 1$ for every even integer $n>4$. This problem remains
still unsolved (for more details, see e.g. [8; Section 2.8 and p.
594]). Another famous conjecture related to prime partitions was
stated by Hardy and Littlewood [5], who predicted the asymptotic
form of $Q_2(n)$ for large even $n$. They conjectured that
\begin{eqnarray}\label{hardlit}
& & Q_2(n)\sim 2C_2\left(\prod_{p_k\in\mathcal{P}, p_k\mid n}
\frac{p_k-1}{p_k-2}\right)\int_2^n\frac{du}{\log^2{u}} \nonumber
\\
 & & \sim 2C_2 \left(\prod_{p_k\in\mathcal{P}, p_k\mid n}
\frac{p_k-1}{p_k-2}\right)\frac{n}{\log^2{n}}, \quad n\to\infty,
\end{eqnarray}
where $C_2$ is the twin prime constant
$$
C_2:=\prod_{p_k\in\mathcal{P}} \left(1-\frac{1}{(p_k-1)^2}\right)
=0.6601618158...
$$
(for the role of $C_2$ in the distribution of the prime numbers,
see again [8; Section 22.20]). This conjecture remains also still
open.

In the present note we do not deal with the asymptotic equivalence
(\ref{hardlit}) but consider the sum function
\begin{equation}\label{summatory}
S(2n)=\sum_{2<k\le n} Q_2(2k), \quad n>2,
\end{equation}
counting all partitions of the even integers from the interval
$(4,2n]$ into two odd prime parts. Sometimes this kind of partitions
are called Goldbach partitions. Let $\Sigma_{2n}$ denote the set
of these partitions. Our main result is the following asymptotic
equivalence.

\begin{theorem}
We have
$$
\mid\Sigma_{2n}\mid=S(2n)\sim\frac{2n^2}{\log^2{n}}, \quad
n\to\infty.
$$
\end{theorem}

Consider now a random experiment. Suppose that we select a
partition uniformly at random from the set $\Sigma_{2n}$, i.e. we
assign the probability $1/S(2n)$ to each Goldbach partition. We
denote by $\mathbb{P}$ the uniform probability measure on
$\Sigma_{2n}$. Let $2X_n\in (4,2n]$ be the number that is
partitioned by this random selection. $2X_n$ is also called the
size of this partition. Using Theorem 1, we determine the limiting
distribution of the random variable $X_n$.

\begin{theorem} If $0<u<1$, then
$$
lim_{n\to\infty}\mathbb{P}\left(\frac{X_n}{n}\le u\right)=u^2.
$$
\end{theorem}

{\it Remark 1.} Using the Prime Number Theorem [8; Section 1.8],
it is easy to show that the number of ordered pairs of primes not
exceeding $2n$ is also $\sim 2n^2/\log^2{n}$; cf. with the result
of Theorem 1. Hence, we conclude that almost all even integers
that are $\le 2n$ have only one partition into two prime parts.

{\it Remark 2.} In probabilistic terms Theorem 2 shows that the
typical size of a random Goldbach partition is a fraction of $2n$.
Moreover, Theorem 2 implies that $X_n/n$ converges weakly, as
$n\to\infty$, to a random variable whose cumulative distribution
function is
$$
F(u) =\left\{\begin{array}{ll} 0 & \qquad  \mbox {if} \qquad u\le
 0, \\
 u^2 & \qquad \mbox {if}\qquad 0<u<1, \\
 1 & \qquad \mbox {if} \qquad u\ge 1.
 \end{array}\right.
$$
 It can be
easily seen that $F(u)$ is the distribution function of
$\max{\{U_1,U_2\}}$, where $U_1$ and $U_2$ are two independent
copies of a uniformly distributed random variable in the interval
$(0,1)$.

{\it Remark 3.} One reason to study the sum function
(\ref{summatory}) is motivated by a result due to Brigham [2]. He
has studied the asymptotic behavior of a similar sum function
related to integer partitions weighted by the sequence of the von
Mangoldt functions (the definition of a von Mangoldt function and
its role in the proof of the Prime Number Theorem may be found in
[8; Section 17.7]). The asymptotic behavior of a single term in
Brigham's sum function was subsequently studied by Richmond [13]
and Yang [17]. Their results are essentially based on Brigham's
observations.

{\it Remark 4.} Another interesting problem on prime partitions is
related to the asymptotic behavior of the coefficients $Q_m(n)$,
the number of prime partitions of $n$ with $m$ parts  (see
(\ref{euler})). Haselgrove and Temperley [9; p. 240] found an
asymptotic form for $Q_m(n)$, whenever $m=m(n)\to\infty$ as
$n\to\infty$ in a proper way. In probabilistic terms their result
can be stated as follows. Consider a random variable, whose
probability distribution function is defined by the ratio
\begin{equation}\label{pdf}
\frac{Q_m(n)}{Q(n)}, \quad m=1,...,n.
\end{equation}
Haselgrove and Temperley [9] showed that this random variable
converges weakly to a non-degenerate random variable as
$n\to\infty$. They also determined the moment generating function
of this limiting variable. The asymptotic form of the mean and the
variance of probability distribution (\ref{pdf}) were found
recently by Ralaivaosaona [12].

Our paper is organized as follows. Section 2 contains some
preliminaries. The proofs of Theorems 1 and 2 are given in Section
3. Our method of proof is essentially based on a classical
Tauberian theorem due to Hardy, Littlewood and Karamata (see [4]).
Finally, in Section 4 we present an extension of our main result.
In particular, we show that the same approach yields similar
results for prime partitions of $n$ into $m>2$ parts whenever $m$
is fixed integer.

\section{Preliminary Results}

We start with a generating function identity for the sequence
$\{Q_2(2k)\}_{k>2}$ of the counts of Goldbach partitions.

\begin{lemma} For any real variable $z$ with $\mid z\mid<1$, let
\begin{equation}\label{ef}
f(z)=\sum_{p_k\in\mathcal{P}} z^{p_k}.
\end{equation}
Then, we have
\begin{equation}\label{ident}
2\sum_{k>2} Q_2(2k)z^{2k} =f^2(z)+f(z^2).
\end{equation}
\end{lemma}

{\it Proof.} Differentiating the left-hand side of (\ref{euler})
twice with respect to $x$ and setting then $x=0$ and $m=2$, we get
\begin{eqnarray}
& & \frac{\partial^2 G(x,z)}{\partial x^2}\mid_{x=0,m=2}
=\sum_{n=1}^\infty z^n\sum_{m=2}^n
m(m-1)Q_m(n)x^{m-2}\mid_{x=0,m=2}
\nonumber \\
& & =2\sum_{n=1}^\infty Q_2(n)z^n =2\sum_{k>2} Q_2(2k)z^{2k}.
\nonumber
\end{eqnarray}
The last equality follows from the obvious identities
$Q_2(1)=Q_2(2)=Q_2(4)=0$ and $Q_2(2k+1)=0$ for $k=1,2,...$. The
right-hand side of (\ref{euler}) can be also written as
$\exp{(-\sum_{p_k\in\mathcal{P}}\log{(1-xz^{p_k})})}$.
Differentiating it twice, in the same way we find that
\begin{eqnarray}
& & \frac{\partial^2 G(x,z)}{\partial x^2}\mid_{x=0}
=\left(exp{\left(-\sum_{p_k\in\mathcal{P}}\log{(1-xz^{p_k})}\right)}
\right) \left(\sum_{p_k\in\mathcal{P}} \frac{z^{p_k}}{1-xz^{p_k}}
\right)^2\mid_{x=0} \nonumber \\
& &
+\left(exp{\left(-\sum_{p_k\in\mathcal{P}}\log{(1-xz^{p_k})}\right)}
\right) \left(\sum_{p_k\in\mathcal{P}}
\frac{z^{2p_k}}{(1-xz^{p_k})^2}\right)\mid_{x=0} \nonumber \\
& & =f^2(z)+f(z^2), \nonumber
\end{eqnarray}
which completes the proof.$\rule{2mm}{2mm}$

Further, we will use a Tauberian theorem by
Hardy-Littlewood-Karamata whose proof may be found in [4; Chapter
7]. We use it in the form given by Odlyzko [11; Section 8.2].

{\bf Hardy-Littlewood-Karamata Theorem.} {\it (See [11; Theorem
8.7, p. 1225].) Suppose that $a_k\ge 0$ for all $k$, and that
$$
g(x)=\sum_{k=0}^\infty a_k x^k
$$
converges for $0\le x<r$. If there is a $\rho>0$ and a function
$L(t)$ that varies slowly at infinity such that
\begin{equation}\label{funcsim}
g(x)\sim (r-x)^{-\rho}L\left(\frac{1}{r-x}\right), \quad x\to r^-,
\end{equation}
then
\begin{equation}\label{sumsim}
\sum_{k=0}^n a_k r^k\sim\left(\frac{n}{r}\right)^\rho
\frac{L(n)}{\Gamma(\rho+1)}, \quad n\to\infty.
\end{equation}}

{\it Remark.} A function $L(t)$ varies slowly at infinity if, for
every $u>0$, $L(ut)\sim L(t)$ as $t\to\infty$.

\section{Proof of the Main Result}

{\it Proof of Theorem 1.} We need to show that power series
(\ref{ident}) satisfies the conditions of
Hardy-Littlewood-Karamata theorem. The next lemma establishes an
asymptotic equivalence of $f(z)$ as $z\to 1^-$.

\begin{lemma} Let $f(z)$ be the power series defined by
(\ref{ef}). Then, as $z\to 1^-$,
$$
f(z)\sim -\frac{1}{\left(\log{\frac{1}{z}}\right)
\left(\log{\log{\frac{1}{z}}}\right)}.
$$
\end{lemma}
{\it Proof.} As usual, by $\pi(y)$ we denote the number of primes
which do not exceed the positive real number $y$. In (\ref{ef}) we
set $z=e^{-t}, t>0,$ and apply an argument similar to that given
by Stong [15] (see also [3]). We have
\begin{eqnarray}\label{intsum}
& & f(e^{-t}) =\int_0^\infty e^{-yt}d\pi(y) =\int_0^\infty
te^{-yt}\pi(y)dy =\int_0^\infty \pi(s/t)e^{-s}ds \nonumber \\
& & =I_1(t)+I_2(t),
\end{eqnarray}
where
$$
I_1(t)=\int_0^{t^{1/2}}\pi(s/t)e^{-s}ds, \quad I_2(t)
=\int_{t^{1/2}}^\infty \pi(s/t)e^{-s}ds.
$$
For $I_1(t)$ we use the bound $\pi(s/t)\le s/t$. Hence, for enough
small $t>0$, we obtain
\begin{eqnarray}\label{ione}
& & 0\le I_1(t) \le\frac{1}{t}\int_0^{t^{1/2}}se^{-s}ds \nonumber \\
& & =\frac{1}{t}\left(-se^{-s}\mid_0^{t^{1/2}}+\int_0^{t^{1/2}}
e^{-s}ds\right) =\frac{1}{t}O(t^{1/2})=O(t^{-1/2}).
\end{eqnarray}
The estimate for $I_2(t)$ follows from the Prime Number Theorem
with an error term given in a suitable form. So, it is known that,
for $y>1$,
$$
\pi(y) =\frac{y}{\log{y}}+O\left(\frac{y}{\log^2{y}}\right)
$$
(see e.g. [10; Theorem 23, p. 65]). Furthermore, for $s\ge
t^{1/2}$, we have $\log{s}\ge -\frac{1}{2}\log{\frac{1}{t}}$.
Hence, as in [15], we get
\begin{eqnarray}\label{pist}
& & \pi(s/t) =\frac{s}{t}\frac{1}{\log{\frac{1}{t}}+\log{s}}
+O\left(\frac{s}{t\left(\log{\frac{1}{t}}+\log{s}\right)^2}\right)
\nonumber \\
& & =\frac{s}{t\log{\frac{1}{t}}}
\left(1+O\left(\frac{\mid\log{s}\mid}{\log{\frac{1}{t}}}\right)\right)
+O\left(\frac{s}{t\log^2{\frac{1}{t}}}\right) \nonumber \\
& & =\frac{s}{t\log{\frac{1}{t}}}
+O\left(\frac{s(1+\mid\log{s}\mid)}{t\log^2{\frac{1}{t}}}\right).
\end{eqnarray}
We also recall that in (\ref{ione}) we have used the obvious
estimate
\begin{equation}\label{estim}
\int_0^{t^{1/2}} se^{-s}ds =O(t^{1/2}).
\end{equation}
Combining (\ref{pist}) and (\ref{estim}), we obtain
\begin{eqnarray}\label{itwo}
& & I_2(t) =\frac{1}{t\log{\frac{1}{t}}} \int_{t^{1/2}}^\infty
se^{-s}ds +O\left(\frac{1}{t\log^2{\frac{1}{t}}}
\int_{t^{1/2}}^\infty s(1+\mid\log{s}\mid)e^{-s}ds\right)
\nonumber \\
& & =\frac{1}{t\log{\frac{1}{t}}}\left(\int_0^\infty se^{-s}ds
+O(t^{1/2})\right) +O\left(\frac{1}{t\log^2{\frac{1}{t}}}\right)
\nonumber \\
& & =\frac{1}{t\log{\frac{1}{t}}}
+O\left(\frac{1}{t^{1/2}\log{\frac{1}{t}}}\right)
+O\left(\frac{1}{t\log^2{\frac{1}{t}}}\right) \nonumber \\
& & \sim\frac{1}{t\log{\frac{1}{t}}}, \quad t\to 0^+.
\end{eqnarray}
Hence, by (\ref{intsum}), (\ref{ione}) and (\ref{itwo}),
$$
f(e^{-t})\sim\frac{1}{t\log{\frac{1}{t}}}, \quad t\to 0^+.
$$
The proof is now completed after the substitution
$t=\log{\frac{1}{z}}$.$\rule{2mm}{2mm}$

Since
$$
\log{\frac{1}{z}}=-\log{z} =-\log{(1-(1-z))} \sim 1-z, \quad z\to
1^-,
$$
the asymptotic equivalence in Lemma 2 becomes
$$
f(z)\sim\frac{1}{(1-z)\log{\frac{1}{1-z}}}, \quad z\to 1^-.
$$
Therefore,
$$
f^2(z)+f(z^2)\sim\frac{1}{(1-z)^2\log^2{\frac{1}{1-z}}}, \quad
z\to 1^-,
$$
which implies that the series $\sum_{k>2} Q(2k)z^{2k}$
satisfies condition (\ref{funcsim}) of Hardy-Littlewood-Karamata
Tauberian theorem with $r=1, \rho=2$ and
$L(t)=\frac{1}{\log^2{t}}$ (see also (\ref{ident})). The
asymptotic equivalence of Theorem 1 follows immediately from
(\ref{sumsim}).$\rule{2mm}{2mm}$

{\it Proof of Theorem 2.} Recall that $2X_n\in (4,2n]$ equals the
size of a Goldbach partition that is chosen uniformly at random
from the set $\Sigma_{2n}$ of all such partitions. Since
$S(2n)=\mid\Sigma_{2n}\mid$ and, for any $N\in (2,n]$,
$S(2N)=\mid\Sigma_{[2N]}\mid$ ($[a]$ denotes the integer part of
the real number $a$), from (\ref{summatory}) it follows that
\begin{equation}\label{probab}
\mathbb{P}(2X_n\le 2N) =\frac{S(2N)}{S(2n)}.
\end{equation}
Setting $N\sim un, 0<u<1,$ and applying Theorem 1 twice - to the
numerator and the denominator of (\ref{probab}), we see that the
limit of (\ref{probab}), as $n\to\infty$, is $u^2$. This completes
the proof.$\rule{2mm}{2mm}$

\section{Prime Partitions with More Than Two Parts}

Let $m>2$ be an integer and let $\Sigma_{m,n}$ denote the set of
prime partitions of the integers from the interval $(4,n]$ into
$m$ parts. The goal of this section is to extend the results of
Theorems 1 and 2 to prime partitions from the class
$\Sigma_{m,n}$. We state them below as two separate theorems.

\begin{theorem} For any fixed integer $m>2$, we have
$$
\mid\Sigma_{m,n}\mid\sim\frac{1}{m!}
\left(\frac{n}{\log{n}}\right)^m, \quad n\to\infty.
$$
\end{theorem}

Furthermore, let $X_{m,n}$ denote the size of a prime partition
selected uniformly at random from the class $\Sigma_{m,n}$. (The
uniform probability measure on $\Sigma_{m,n}$ is again denoted by
$\mathbb{P}$.)

\begin{theorem} If $0<u<1$ and $m$ is as in Theorem 3, then
$$
\lim_{n\to\infty}\mathbb{P}\left(\frac{X_{m,n}}{n}\le u\right)
=u^m.
$$
\end{theorem}

 Theorem 4 shows a weak convergence similar to that established in
 Theorem 2. Namely, for any fixed integer $m$, $X_{m,n}/n$ converges,
 as $n\to\infty$, to $\max{\{U_1,...,U_m\}}$, where $U_1,...,U_m$
 are independent copies of a random variable that is uniformly
 distributed in the interval $(0,1)$.

 Below we only sketch the proof of Theorem 3. The proof of Theorem
 4 is almost identical with that of Theorem 2.

 {\it Sketch of the proof of Theorem 3.} Our main tool is again the
 generating function identity (\ref{euler}). We notice first that
 the coefficients $Q_m(n)$ are $=0$ if either $m$ is odd and $n$
 is even or $m$ is even and $n$ is odd. By the definition of
 $Q_m(n)$, we also have
 \begin{equation}\label{summatorym}
 \mid\Sigma_{m,n}\mid =\sum_{k\le n} Q_m(k).
 \end{equation}
 We compute the $m$th derivative of the infinite product in (\ref{euler})
 using Faa di Bruno formula for derivatives of compound
 functions (see e.g. [14; Section 2.8]). We introduce the following
 auxiliary notations:
 \begin{eqnarray}\label{ab}
 & & b(x)=b(x,z):=-\sum_{p_k\in\mathcal{P}}
 \log{(1-xz^{p_k})}, \nonumber \\
 & & b_j=b_j(x,z):=\frac{\partial^j b(x,z)}{\partial x^j},
 j=1,...,m.
 \end{eqnarray}
 Using formulae (43) and (46) of [14; Section 2.8], we obtain
 \begin{equation}\label{mthab}
 \frac{d^m}{dx^m}e^{b(x)} =e^{b(x)} b_1^m+R_m,
 \end{equation}
 where
 \begin{equation}\label{bruno}
 R_m=R_m(x,z) =e^{b(x,z)}\widetilde{\sum} \frac{m!}{k_1!...k_m!}
 \left(\frac{b_1}{1!}\right)^{k_1}...
 \left(\frac{b_m}{m!}\right)^{k_m}
 \end{equation}
 and  $\widetilde{\sum}$ denotes the sum over all integers $k_j\ge
 0, j=1,...,m$, such that $\sum_{j=1}^m jk_j=m$ and $k_1<m$.
 Setting $x=0$ in (\ref{ab}), we find that $b(0,z)=0$ and $b_j(0,z)=f(z^j),
 j=1,...,m,$ where the function $f(z)$ is defined (\ref{ef}).
 Moreover, in the right-hand side of (\ref{bruno}) we have
 $k_1+...+k_m\le m-1$. In fact, since $k_1<m$ by the definition of
 $\widetilde{\sum}$ at least one $k_j$ is $>0$ for $j\ge 2$. Hence
 if $m=k_1+k_2+...+k_m$, then $m<k_1+\sum_{j=2}^m jk_j=m$.
 Since $f(z^j)=O(f(z))$ as $z\to 1^-$ and since
 $k_1+...+k_m\le m-1$, we conclude that $R_m(0,z)=O(f^{m-1}(z))$. Therefore
 (\ref{mthab}) becomes
 $$
 \frac{d^m}{dx^m} e^{b(x)}=f^m(z)+O(f^{m-1}(z)),
 $$
 or, equivalently,
 \begin{equation}\label{mth}
 \frac{\partial^m G(x,z)}{\partial x^m}\mid_{x=0} =f^m(z)+O(f^{m-1}(z))
 \end{equation}
 as $z\to 1^-$. On the other hand,
 \begin{equation}\label{mthtwo}
\frac{\partial^m G(x,z)}{\partial x^m}\mid_{x=0} =m!\sum_{k\ge m}
Q_m(k)z^k.
\end{equation}
Applying Lemma 2, as in the proof of Theorem 1, we obtain the
asymptotic equivalence
\begin{equation}\label{mthasymp}
f^m(z) \sim\frac{1}{(1-z)^m\log^m{\frac{1}{1-z}}}, \quad z\to 1^-.
\end{equation}
The observations in (\ref{mth})-(\ref{mthasymp}) imply that
$$
\frac{\partial^m G(x,z)}{\partial x^m}\mid_{x=0}
\sim\frac{1}{(1-z)^m\log^m{\frac{1}{1-z}}}, \quad z\to 1^-.
$$
So, condition (\ref{funcsim}) of Hardy-Littlewood-Karamata theorem
is satisfied with $r=1, \rho=m$ and $L(t)=\frac{1}{\log^m{t}}$.
The required result follows at once from (\ref{summatorym}) and
(\ref{sumsim}).$\rule{2mm}{2mm}$

\section*{Acknowledgements}

I am grateful to the referee for carefully reading the paper and
for his helpful comments and suggestions.

\end{document}